\documentclass[12pt]{article}

\usepackage{amsmath, epsfig, cite}
\usepackage{amssymb}
\usepackage{amsfonts}
\usepackage{latexsym}
\usepackage{amsthm}

\newtheorem{thm}{Theorem}[section]

\numberwithin{equation}{section}

\renewcommand{\thefootnote}

\begin{document}

\begin{center}
{\large\bf On a conjectural series of Sun for the mathematical\\
constant $\beta(4)$
 \footnote{ The work is supported by the National Natural Science Foundation of China (No. 12071103).}}
\end{center}

\renewcommand{\thefootnote}{$\dagger$}

\vskip 2mm \centerline{Chuanan Wei}
\begin{center}
{School of Biomedical Information and Engineering\\ Hainan Medical
University, Haikou 571199, China
\\
 Email address: weichuanan78@163.com}
\end{center}


\vskip 0.7cm \noindent{\bf Abstract.} Series expansions for the
mathematical constant $\beta(4)$ are rare in the history. With the
help of the operator method and a hypergeometric transformation, we
prove a surprising conjectural series of Sun for $\beta(4)$.
Furthermore, we find five new series for the same constant in this
paper.



\section{Introduction}

For a nonnegative integer $n$, define the shifted-factorial to be
\begin{align*}
(x)_0=1 \quad\text{and}\quad (x)_n=x(x+1)\cdots(x+n-1),
\end{align*}
where $x$ is a complex variable. For two complex sequences
$\{a_i\}_{i\geq1}$ and $\{b_j\}_{j\geq1}$, define the hypergeometric
series by
$$
_{r+1}F_{r}\left[\begin{array}{c}
a_1,a_2,\ldots,a_{r+1}\\
b_1,b_2,\ldots,b_{r}
\end{array};\, z
\right] =\sum_{k=0}^{\infty}\frac{(a_1)_k(a_2)_k\cdots(a_{r})_k}
{(1)_k(b_1)_k\cdots(b_{r})_k}z^k,
$$
where $z$ is any complex number. Then Whipple's transformation from
a $_7F_6$ series to a $_4F_3$ series(cf.\cite[P. 28]{Bailey}) can be
expressed as
\begin{align}
&{_{7}F_{6}}\left[\begin{array}{cccccccc}
  a,1+\frac{a}{2},b,c,d,e,-n\\
  \frac{a}{2},1+a-b,1+a-c,1+a-d,1+a-e,1+a+n
\end{array};1
\right]
\notag\\
&=\frac{(1+a)_n(1+a-d-e)_n}{(1+a-d)_n(1+a-e)_n}
{_{4}F_{3}}\left[\begin{array}{cccccccc}
  1+a-b-c,d,e,-n\\
  1+a-b,1+a-c,d+e-a-n
\end{array};1
\right]. \label{Whipple}
\end{align}

 Recall two series for $1/\pi$:
\begin{align}
&\sum_{k=0}^{\infty}\frac{1+6k}{4^k}\frac{(\frac{1}{2})_k^3}{(1)_k^3}
=\frac{4}{\pi}, \label{Ramanujan}
\end{align}
\begin{align}
\label{Gosper}
&\sum_{k=0}^{\infty}\frac{1+6k}{(-8)^k}\frac{(\frac{1}{2})_k^3}{(1)_k^3}
=\frac{2\sqrt{3}}{\pi}.
\end{align}
Here \eqref{Ramanujan} is due to Ramanujan \cite{Ramanujan} and
\eqref{Gosper} is a special case of Gosper's $_4F_3$ summation
 (cf. \cite{Gosper}):
$$
_{4}F_{3}\left[\begin{array}{c}
a,1+\frac{a}{3},b,1-b\\[2mm]
\frac{a}{3},\frac{1+a+b}{2},1+\frac{a-b}{2}
\end{array};\, -\frac{1}{8}
\right]
=\frac{2^a\Gamma(\frac{1+a+b}{2})\Gamma(1+\frac{a-b}{2})}{\Gamma(\frac{1}{2})\Gamma(1+a)},
$$
where $\Gamma(x)$ is the famous Gamma function. Let $G$ be the
Catalan constant
$$G=\beta(2)=\sum_{k=0}^{\infty}\frac{(-1)^k}{(1+2k)^2},$$
where $\beta(z)$ is the Dirichlet beta function given by
$$\beta(z)=\sum_{k=0}^{\infty}\frac{(-1)^k}{(1+2k)^z}.$$
 Two interesting results due to Guillera \cite{Guillera-b} can be
 stated as
\begin{align}
&\sum_{k=0}^{\infty}\frac{2+3k}{(-8)^k}\frac{(1)_k^3}{(\frac{3}{2})_k^3}
=2G, \label{Guillera-a}\\\label{Guillera-b}
&\sum_{k=0}^{\infty}\frac{3+4k}{(-1)^k}\frac{(1)_k^3}{(\frac{3}{2})_k^3}
=2G.
\end{align}

In 2021, Guo and Lian \cite{Guo-c} conjectured two double series for
$\pi$ related to \eqref{Ramanujan} and \eqref{Gosper}:
\begin{align}
&\:\:\:\sum_{k=0}^{\infty}\frac{1+6k}{4^k}\frac{(\frac{1}{2})_k^3}{(1)_k^3}
\sum_{j=1}^{k}\bigg\{\frac{1}{(2j-1)^2}-\frac{1}{16j^2}\bigg\}
=\frac{\pi}{12},
 \label{guo-a}\\\label{guo-b}
&\sum_{k=0}^{\infty}\frac{1+6k}{(-8)^k}\frac{(\frac{1}{2})_k^3}{(1)_k^3}
\sum_{j=1}^{k}\bigg\{\frac{1}{(2j-1)^2}-\frac{1}{16j^2}\bigg\}
=-\frac{\sqrt{2}}{48}\pi,
\end{align}
which have been proved by the author \cite{Wei-a}. For two positive
integers $\ell$ and $n$, define the generalized harmonic numbers of
order $\ell$ as
\[H_{n}^{(\ell)}(x)=\sum_{k=1}^n\frac{1}{(x+k)^{\ell}},\]
where  $x$ is a complex variable. The $x=0$ case of them become the
harmonic numbers of order $\ell$:
\[H_{n}^{(\ell)}=\sum_{k=1}^n\frac{1}{k^{\ell}}.\]
Recently,  Sun \cite{Sun-c} rewrote \eqref{guo-a} and \eqref{guo-b}
as
\begin{align*}
&\:\:\:\sum_{k=0}^{\infty}\frac{1+6k}{4^k}\frac{(\frac{1}{2})_k^3}{(1)_k^3}
\bigg\{H_{2k}^{(2)}-\frac{5}{16}H_{k}^{(2)}\bigg\}
=\frac{\pi}{12},\\[2mm]
&\sum_{k=0}^{\infty}\frac{1+6k}{(-8)^k}\frac{(\frac{1}{2})_k^3}{(1)_k^3}\bigg\{H_{2k}^{(2)}-\frac{5}{16}H_{k}^{(2)}\bigg\}
=-\frac{\sqrt{2}\,\pi}{48}.
\end{align*}
In the same paper, he proposed the following two conjectures
associated with \eqref{Guillera-a} and \eqref{Guillera-b} (cf.
\cite[Equations (3.24) and (3.37)]{Sun-c}) :
\begin{align}
&\sum_{k=0}^{\infty}\frac{2+3k}{(-8)^k}\frac{(1)_k^3
}{(\frac{3}{2})_k^3
}\bigg\{H_{1+2k}^{(2)}-\frac{5}{4}H_{k}^{(2)}\bigg\} =2\beta(4),
\label{Sun-a}\\[2mm]
&\sum_{k=0}^{\infty}\frac{3+4k}{(-1)^k}\frac{(1)_k^3}{(\frac{3}{2})_k^3}
\bigg\{H_{1+2k}^{(2)}-\frac{1}{2}H_{k}^{(2)}\bigg\} =2\beta(4),
\label{Sun-b}
\end{align}
where $\beta(4)$ denotes the mathematical constant
$$\beta(4)=\sum_{k=0}^{\infty}\frac{(-1)^k}{(1+2k)^4}.$$
Equation \eqref{Sun-b} has been verified by Au \cite{Au}, but the
surprising series \eqref{Sun-a} is still open. For more series on
mathematical constants, we refer the reader to the papers
\cite{Schlosser,Sun-a,Wang,Zudilin}. Inspired by the works just
mentioned, we shall established the following theorem.

\begin{thm}\label{thm-a}
Equation \eqref{Sun-a} is true.
\end{thm}

Furthermore, we shall dispaly the following five new series for the
mathematical constant $\beta(4)$.
\begin{thm}\label{thm-b}
\begin{align}
\sum_{k=0}^{\infty}\frac{(1)_k}{(\frac{3}{2})_k(1+2k)}\bigg\{H_{1+2k}^{(2)}-\frac{1}{4}H_{k}^{(2)}\bigg\}
=2\beta(4).
 \label{eq:wei-b}
\end{align}
\end{thm}

\begin{thm}\label{thm-c}
\begin{align}
\sum_{k=0}^{\infty}\bigg(\frac{1}{4}\bigg)^k\frac{(1)_k^2}{(\frac{5}{4})_k(\frac{7}{4})_k}\bigg\{\frac{5+6k}{1+2k}H_{k}^{(2)}-\frac{24+32k}{(1+2k)^3}\bigg\}
=-24\beta(4).
 \label{eq:wei-c}
\end{align}
\end{thm}

\begin{thm}\label{thm-d}
\begin{align}
&\sum_{k=0}^{\infty}\bigg(\frac{-1}{4}\bigg)^k\frac{(1)_k^3(\frac{1}{2})_k}{(\frac{5}{4})_k^2(\frac{7}{4})_k^2}
\Big\{(19+56k+40k^2)\big[H_{k}^{(2)}+H_{1+2k}^{(2)}-4H_{3+4k}^{(2)}\big]+4\Big\}
\notag\\
&\:\:=-72\beta(4).
 \label{eq:wei-d}
\end{align}
\end{thm}

\begin{thm}\label{thm-e}
\begin{align}
\sum_{k=0}^{\infty}\bigg(\frac{16}{27}\bigg)^k\frac{(1)_k^2}{(\frac{7}{6})_k(\frac{11}{6})_k}
\bigg\{\frac{21+22k}{1+2k}H_{1+2k}^{(2)}-\frac{145+174k}{(1+2k)^3}\bigg\}
=-120\beta(4).
 \label{eq:wei-e}
\end{align}
\end{thm}

\begin{thm}\label{thm-f}
\begin{align}
&\sum_{k=0}^{\infty}\bigg(\frac{-1}{27}\bigg)^k\frac{(1)_k^2}{(\frac{7}{6})_k(\frac{11}{6})_k}
\bigg\{\frac{83+192k+112k^2}{64(1+4k)(3+4k)}H_{k}^{(2)}-\frac{10+39k+48k^2+16k^3}{(1+4k)^3(3+4k)^3(5+6k)^{-1}}\bigg\}
\notag\\
&\:\:=-\frac{15}{8}\beta(4).
 \label{eq:wei-f}
\end{align}
\end{thm}

According to the limit of the difference operator, we may define the
derivative operator $\mathcal{D}_{x}$ by
\begin{align*}
\mathcal{D}_{x}f(x)=\lim_{h\to0}\frac{f(x+h)-f(x)}{h}.
\end{align*}
Then there are two relations:
\begin{align*}
&\mathcal{D}_{x}(x)_n=(x)_nH_n(x-1),\\
&\:\:\mathcal{D}_{x}H_n(x)=-H_n^{(2)}(x).
\end{align*}

The rest of the paper is arranged as follows. We shall prove Theorem
\ref{thm-a} via the derivative operator and a hypergeometric
transformation in Section 2. We shall certify Theorem \ref{thm-b}
through the derivative operator and \eqref{Whipple} in Section 3.
Similarly, The proof of Theorems \ref{thm-c}-\ref{thm-f} will be
provided in Section 4.

\section{Proof of Theorem \ref{thm-a}}
For proving Theorem \ref{thm-a}, we draw support on
 the following hypergeometric  transformation (cf. \cite[Theorem
27]{Chu-b}):
\begin{align}
&\sum_{k=0}^{\infty}\frac{(b)_k(c)_k(d)_k(1+a-b-c)_k(1+a-b-d)_{k}(1+a-c-d)_{k}}{(b+e-a)_{k}(c+e-a)_{k}(d+e-a)_{k}}
\notag\\[1mm]
&\quad\times\frac{(e)_{3k}}{(1+a-b)_{2k}(1+a-c)_{2k}(1+a-d)_{2k}}\mu_k(a,b,c,d,e)
\notag\\[1mm]
&\:+\frac{\Gamma(1+a-b)\Gamma(1+a-c)\Gamma(1+a-d)\Gamma(1+a-e)}{\Gamma(b)\Gamma(c)\Gamma(d)\Gamma(e)}
\notag\\[1mm]
&\quad\times\frac{\Gamma(b+e-a)\Gamma(c+e-a)\Gamma(d+e-a)\Gamma(1+2a-b-c-d-e)}{\Gamma(1+a-b-c)\Gamma(1+a-b-d)\Gamma(1+a-c-d)}
\notag\\[1mm]
&\:=\sum_{k=0}^{\infty}(a+2k)\frac{(b)_k(c)_k(d)_k(e)_k}{(1+a-b)_{k}(1+a-c)_{k}(1+a-d)_{k}(1+a-e)_{k}},
\label{equation-a}
\end{align}
where
\begin{align*}
\mu_k(a,b,c,d,e)&=\frac{(a-c+2k)(a-e)}{a-c-e-k}
-\frac{(c+k)(e+3k)(a-e)(1+a-b-d+k)}{(1+a-d+2k)(a-b-e-k)(a-c-e-k)}
\\[1mm]
&\quad+\frac{(b+k)(c+k)(e+3k)(1+e+3k)}{(1+a-b+2k)(1+a-c+2k)(1+a-d+2k)}
\\[1mm]
&\qquad\times\frac{(a-e)(1+a-b-d+k)(1+a-c-d+k)}{(a-b-e-k)(a-c-e-k)(a-d-e-k)}.
\end{align*}

\begin{proof}[{\bf{Proof of Theorem \ref{thm-a}}}]

Choose $(a,b,c,d,e)=(1,x,1-x,\frac{1}{2},-n)$ in \eqref{equation-a}
to obtain
\begin{align}
&\sum_{k=0}^{n}\frac{(x)_k(1-x)_k(\frac{1}{2}+x)_k(\frac{3}{2}-x)_k(\frac{1}{2})_k(1)_k}{(1+x)_{2k}(2-x)_{2k}(\frac{3}{2})_{2k}}
\frac{(-n)_{3k}}{(x-1-n)_k(-x-n)_{k}(-\frac{1}{2}-n)_{k}}
\notag\\[1mm]
&\quad\times
A_k(x;n)=\sum_{k=0}^{n}\frac{(x)_k(1-x)_k(-n)_k}{(2-x)_{k}(1+x)_{k}(2+n)_{k}},
\label{eq:wei-aa}
\end{align}
where
\begin{align*}
A_k(x;n)
&=\frac{(x+2k)(1+n)}{x-k+n}-\frac{(1-x+k)(3-2x+2k)(3k-n)(1+n)}{(3+4k)(x-k+n)(1-x-k+n)}
\\[2mm]
&+\frac{(x+k)(1-x+k)(1+2x+2k)(3-2x+2k)(3k-n)(1+3k-n)(1+n)}{(1+x+2k)(2-x+2k)(3+4k)(1-2k+2n)(x-k+n)(1-x-k+n)}.
\end{align*}
Apply the operator $\mathcal{D}_{x}$ on \eqref{eq:wei-aa} to get
\begin{align}
&\sum_{k=0}^{n}\frac{(x)_k(1-x)_k(\frac{1}{2}+x)_k(\frac{3}{2}-x)_k(\frac{1}{2})_k(1)_k}{(1+x)_{2k}(2-x)_{2k}(\frac{3}{2})_{2k}}
\frac{(-n)_{3k}}{(x-1-n)_k(-x-n)_{k}(-\frac{1}{2}-n)_{k}}
\notag\\[2mm]
&\quad\times\Big\{H_{k}(x-1)-H_{k}(-x)+H_{k}(x-\tfrac{1}{2})-H_{k}(\tfrac{1}{2}-x)+H_{2k}(1-x)-H_{2k}(x)
\notag\\[1mm]
&\qquad+H_{k}(-1-x-n)-H_{k}(x-2-n)\Big\}A_k(x;n)
\notag\\[1mm]
&\:+\sum_{k=0}^{n}\frac{(x)_k(1-x)_k(\frac{1}{2}+x)_k(\frac{3}{2}-x)_k(\frac{1}{2})_k(1)_k}{(1+x)_{2k}(2-x)_{2k}(\frac{3}{2})_{2k}}
\frac{(-n)_{3k}\mathcal{D}_{x}A_k(x;n)}{(x-1-n)_k(-x-n)_{k}(-\frac{1}{2}-n)_{k}}
\notag\\[1mm]
&\:\:=
\sum_{k=0}^{n}\frac{(x)_k(1-x)_k(-n)_k}{(2-x)_{k}(1+x)_{k}(2+n)_{k}}
\Big\{H_{k}(x-1)-H_{k}(-x)+H_{k}(1-x)-H_{k}(x)\Big\}.
\label{eq:wei-bb}
\end{align}
Dividing both sides by $1-2x$ and noticing the relation
 \begin{align*}
 \frac{1}{v-u-2x}\Big\{H_m(x+u)-H_m(v-x)\Big\}=\sum_{i=1}^m\frac{1}{(x+u+i)(v-x+i)},
\end{align*}
which will frequently be used elsewhere without indication,
\eqref{eq:wei-bb} can be restated as
\begin{align*}
&\sum_{k=0}^{n}\frac{(x)_k(1-x)_k(\frac{1}{2}+x)_k(\frac{3}{2}-x)_k(\frac{1}{2})_k(1)_k}{(1+x)_{2k}(2-x)_{2k}(\frac{3}{2})_{2k}}
\frac{(-n)_{3k}}{(x-1-n)_k(-x-n)_{k}(-\frac{1}{2}-n)_{k}}
\notag\\[1mm]
&\quad\times\bigg\{\sum_{i=1}^k\frac{1}{(x-1+i)(-x+i)}+\sum_{i=1}^k\frac{1}{(x-\frac{1}{2}+i)(\frac{1}{2}-x+i)}-\sum_{i=1}^{2k}\frac{1}{(x+i)(1-x+i)}
\notag\\[1mm]
&\qquad-\sum_{i=1}^k\frac{1}{(x-2-n+i)(-1-x-n+i)}\bigg\}A_k(x;n)
\notag\\[1mm]
&\:+\sum_{k=0}^{n}\frac{(x)_k(1-x)_k(\frac{1}{2}+x)_k(\frac{3}{2}-x)_k(\frac{1}{2})_k(1)_k}{(1+x)_{2k}(2-x)_{2k}(\frac{3}{2})_{2k}}
\frac{(-n)_{3k}\mathcal{D}_{x}A_k(x;n)}{(x-1-n)_k(-x-n)_{k}(-\frac{1}{2}-n)_{k}(1-2x)}
\notag\\[1mm]
&\:\:=
\sum_{k=0}^{n}\frac{(x)_k(1-x)_k(-n)_k}{(2-x)_{k}(1+x)_{k}(2+n)_{k}}
\bigg\{\sum_{i=1}^k\frac{1}{(x-1+i)(-x+i)}-\sum_{i=1}^k\frac{1}{(x+i)(1-x+i)}\bigg\}.
\end{align*}
Taking $(x,n)\to(\frac{1}{2},\infty)$ in the last equation, it is
not difficult to see that
\begin{align}
&\sum_{k=0}^{\infty}\frac{(1)_{2k}^3}{(\frac{3}{2})_{2k}^38^{2k}}
\bigg\{\frac{49+342k+840k^2+880k^3+336k^4}{(3+4k)^3}\Big[H_{1+4k}^{(2)}-\frac{5}{4}H_{2k}^{(2)}-1\Big]
\notag\\
&\qquad+\frac{4(1+k)(1+2k)(2+3k)(5+6k)}{(3+4k)^3}\bigg\}=2\beta(4)-2G.
 \label{eq:wei-cc}
\end{align}
The $(x,n)\to(\frac{1}{2},\infty)$ case of \eqref{eq:wei-aa} reads
\begin{align}
\sum_{k=0}^{\infty}\frac{(1)_{2k}^3}{(\frac{3}{2})_{2k}^38^{2k}}
\frac{49+342k+840k^2+880k^3+336k^4}{(3+4k)^3}=2G.
 \label{eq:wei-dd}
\end{align}
The sum of \eqref{eq:wei-cc} and \eqref{eq:wei-dd} produces
\begin{align}
&\sum_{k=0}^{\infty}\frac{(1)_{2k}^3}{(\frac{3}{2})_{2k}^38^{2k}}
\bigg\{\frac{49+342k+840k^2+880k^3+336k^4}{(3+4k)^3}\Big[H_{1+4k}^{(2)}-\frac{5}{4}H_{2k}^{(2)}\Big]
\notag\\
&\qquad+\frac{4(1+k)(1+2k)(2+3k)(5+6k)}{(3+4k)^3}\bigg\}=2\beta(4).
 \label{eq:wei-ee}
\end{align}
Be means of the parity, \eqref{eq:wei-ee} can be manipulated as
\begin{align*}
&\sum_{k=0}^{\infty}\frac{2+6k}{8^{2k}}\frac{(1)_{2k}^3}{(\frac{3}{2})_{2k}^3}\Big[H_{1+4k}^{(2)}-\frac{5}{4}H_{2k}^{(2)}\Big]
-\sum_{k=0}^{\infty}\frac{5+6k}{8^{1+2k}}\frac{(1)_{1+2k}^3}{(\frac{3}{2})_{1+2k}^3}\Big[H_{3+4k}^{(2)}-\frac{5}{4}H_{1+2k}^{(2)}\Big]
\\
&\:=\beta(4).
\end{align*}
Therefore, we complete the proof of Theorem \ref{thm-a}.
\end{proof}

\section{Proof of Theorem \ref{thm-b}}

Now we are ready to prove Theorem \ref{thm-b}.

\begin{proof}[{\bf{Proof of Theorem \ref{thm-b}}}]
Employ the derivative operator $\mathcal{D}_{b}$ on \eqref{Whipple}
to gain
\begin{align}
&\sum_{k=0}^{n}\frac{(a)_k(1+\frac{a}{2})_k(b)_k(c)_k(d)_k(e)_k(-n)_k}{(1)_k(\frac{a}{2})_k(1+a-b)_k(1+a-c)_k(1+a-d)_k(1+a-e)_k(1+a+n)_k}
\notag\\[2pt]
 &\quad\times\:
 \Big\{H_k(b-1)+H_k(a-b)\Big\}
 \notag\\[2pt]
&\:=\frac{(1+a)_n(1+a-d-e)_n}{(1+a-d)_n(1+a-e)_n}
\sum_{k=0}^{n}\frac{(1+a-b-c)_k(d)_k(e)_k(-n)_k}{(1)_k(1+a-b)_k(1+a-c)_k(d+e-a-n)_k}
\notag\\[2pt]
  &\quad\times\:\Big\{H_k(a-b)-H_k(a-b-c)\Big\}.
\label{eq:wei-ff}
\end{align}
Applying the derivative operator $\mathcal{D}_{b}$ on
\eqref{eq:wei-ff} and then fixing $c=b$, it is easy to show that
\begin{align}
&\sum_{k=0}^{n}\frac{(a)_k(1+\frac{a}{2})_k(b)_k^2(d)_k(e)_k(-n)_k}{(1)_k(\frac{a}{2})_k(1+a-b)_k^2(1+a-d)_k(1+a-e)_k(1+a+n)_k}
\notag\\[2pt]
 &\quad\times\:
 \Big\{\big[H_k(b-1)+H_k(a-b)\big]^2-\big[H_k^{(2)}(b-1)-H_k^{(2)}(a-b)\big]\Big\}
 \notag\\[2pt]
&\:=\frac{(1+a)_n(1+a-d-e)_n}{(1+a-d)_n(1+a-e)_n}
\sum_{k=0}^{n}\frac{(1+a-2b)_k(d)_k(e)_k(-n)_k}{(1)_k(1+a-b)_k^2(d+e-a-n)_k}
\notag\\[2pt]
  &\quad\times\:
  \Big\{\big[H_k(a-b)-H_k(a-2b)\big]^2+\big[H_k^{(2)}(a-b)-H_k^{(2)}(a-2b)\big]\Big\}.
\label{eq:wei-gg}
\end{align}
Applying the derivative operator $\mathcal{D}_{c}$ on
\eqref{eq:wei-ff} and then setting $c=b$, it is routine to discover
that
\begin{align}
&\sum_{k=0}^{n}\frac{(a)_k(1+\frac{a}{2})_k(b)_k^2(d)_k(e)_k(-n)_k}{(1)_k(\frac{a}{2})_k(1+a-b)_k^2(1+a-d)_k(1+a-e)_k(1+a+n)_k}
\notag\\[2pt]
 &\quad\times\:
 \Big\{H_k(b-1)+H_k(a-b)\Big\}^2
 \notag\\[2pt]
&\:=\frac{(1+a)_n(1+a-d-e)_n}{(1+a-d)_n(1+a-e)_n}
\sum_{k=0}^{n}\frac{(1+a-2b)_k(d)_k(e)_k(-n)_k}{(1)_k(1+a-b)_k^2(d+e-a-n)_k}
\notag\\[2pt]
  &\quad\times\:
  \Big\{\big[H_k(a-b)-H_k(a-2b)\big]^2-H_k^{(2)}(a-2b)\Big\}.
\label{eq:wei-hh}
\end{align}
The difference of \eqref{eq:wei-gg} and \eqref{eq:wei-hh} engenders
\begin{align*}
&\sum_{k=0}^{n}\frac{(a)_k(1+\frac{a}{2})_k(b)_k^2(d)_k(e)_k(-n)_k}{(1)_k(\frac{a}{2})_k(1+a-b)_k^2(1+a-d)_k(1+a-e)_k(1+a+n)_k}
 \end{align*}
\begin{align*}
 &\quad\times\:
 \Big\{H_k^{(2)}(a-b)-H_k^{(2)}(b-1)\Big\}
\notag\\[2pt]
&\:=\frac{(1+a)_n(1+a-d-e)_n}{(1+a-d)_n(1+a-e)_n}
\sum_{k=0}^{n}\frac{(1+a-2b)_k(d)_k(e)_k(-n)_k}{(1)_k(1+a-b)_k^2(d+e-a-n)_k}
  H_k^{(2)}(a-b).
\end{align*}
The $(a,b,d,e)=(\frac{3}{2},1,1,1)$ case of this identity is
\begin{align}
&\sum_{k=0}^{n}(3+4k)\frac{(1)_k^3(-n)_k}{(\frac{3}{2})_k^3(\frac{5}{2}+n)_k}
 \Big\{4H_{1+2k}^{(2)}-2H_k^{(2)}-4\Big\}
 \notag\\[2pt]
&\:=\frac{3+2n}{1+2n}\sum_{k=0}^{n}\frac{(\frac{1}{2})_k(1)_k(-n)_k}{(\frac{3}{2})_k^2(\frac{1}{2}-n)_k}
 \Big\{4H_{1+2k}^{(2)}-H_k^{(2)}-4\Big\}.
\label{eq:wei-ii}
\end{align}
In terms of \eqref{Whipple}, form \eqref{eq:wei-ii} we may deduce
\begin{align}
&\sum_{k=0}^{n}(3+4k)\frac{(1)_k^3(-n)_k}{(\frac{3}{2})_k^3(\frac{5}{2}+n)_k}
 \Big\{4H_{1+2k}^{(2)}-2H_k^{(2)}\Big\}
 \notag\\[2pt]
&\:=\frac{3+2n}{1+2n}\sum_{k=0}^{n}\frac{(\frac{1}{2})_k(1)_k(-n)_k}{(\frac{3}{2})_k^2(\frac{1}{2}-n)_k}
 \Big\{4H_{1+2k}^{(2)}-H_k^{(2)}\Big\}.
\label{eq:wei-jj}
\end{align}
Taking $n\to\infty$ in \eqref{eq:wei-jj} and utilizing
\eqref{Sun-b}, we arrive at \eqref{eq:wei-b}.
\end{proof}

\section{Proof of Theorems \ref{thm-c}-\ref{thm-f}}

Firstly, we begin to prove Theorem \ref{thm-c}.

\begin{proof}[{\bf{Proof of Theorem \ref{thm-c}}}]
Recall the known hypergeometric transformation (cf. \cite[Theorem
9]{Chu-b}):
\begin{align}
&\sum_{k=0}^{\infty}\frac{(c)_k(d)_k(e)_k(1+a-b-c)_k(1+a-b-d)_{k}(1+a-b-e)_{k}}{(1+a-c)_{k}(1+a-d)_{k}(1+a-e)_{k}(1+2a-b-c-d-e)_{k}}
\notag\\[1mm]
&\quad\times\frac{(-1)^k}{(1+a-b)_{2k}}\nu_k(a,b,c,d,e)
\notag\\[1mm]
&\:=\sum_{k=0}^{\infty}(a+2k)\frac{(b)_k(c)_k(d)_k(e)_k}{(1+a-b)_{k}(1+a-c)_{k}(1+a-d)_{k}(1+a-e)_{k}},
\label{equation-b}
\end{align}
where
\begin{align*}
\nu_k(a,b,c,d,e)&=\frac{(1+2a-b-c-d+2k)(a-e+k)}{1+2a-b-c-d-e+k}
\\[1mm]
&\quad+\frac{(1+a-b-c+k)(1+a-b-d+k)(e+k)}{(1+a-b+2k)(1+2a-b-c-d-e+k)}.
\end{align*}
Choosing $(a,b,c,d,e)=(1,\frac{1}{2},x,1-x,-n)$ in
\eqref{equation-b}, there holds
\begin{align}
&\sum_{k=0}^{n}\bigg(\frac{-1}{4}\bigg)^k\frac{(x)_k(1-x)_k(\frac{1}{2}+x)_k(\frac{3}{2}-x)_k}{(1+x)_{k}(2-x)_{k}(\frac{3}{4})_{k}(\frac{5}{4})_{k}}
\frac{(-n)_k}{(2+n)_{k}}B_k(x;n)
\notag\\[1mm]
&\:\:
=\sum_{k=0}^{n}\frac{(x)_k(1-x)_k(-n)_k}{(2-x)_{k}(1+x)_{k}(2+n)_{k}},
\label{eq:wei-at}
\end{align}
where
\begin{align*}
B_k(x;n)
&=\frac{(3+4k)(1+k+n)}{3+2k+2n}+\frac{(1+2x+2k)(3-2x+2k)(k-n)}{(3+4k)(3+2k+2n)}.
\end{align*}
Employ the operator $\mathcal{D}_{x}$ on \eqref{eq:wei-at} to obtain
\begin{align*}
&\sum_{k=0}^{n}\bigg(\frac{-1}{4}\bigg)^k\frac{(x)_k(1-x)_k(\frac{1}{2}+x)_k(\frac{3}{2}-x)_k}{(1+x)_{k}(2-x)_{k}(\frac{3}{4})_{k}(\frac{5}{4})_{k}}
\frac{(-n)_k}{(2+n)_{k}}B_k(x;n)
\notag\\[1mm]
&\quad\times\Big\{H_{k}(x-1)-H_{k}(-x)+H_{k}(x-\tfrac{1}{2})-H_{k}(\tfrac{1}{2}-x)+H_{k}(1-x)-H_{k}(x)\Big\}
\notag\\[1mm]
&\:+\sum_{k=0}^{n}\bigg(\frac{-1}{4}\bigg)^k\frac{(x)_k(1-x)_k(\frac{1}{2}+x)_k(\frac{3}{2}-x)_k}{(1+x)_{k}(2-x)_{k}(\frac{3}{4})_{k}(\frac{5}{4})_{k}}
\frac{(-n)_k}{(2+n)_{k}} \mathcal{D}_{x}B_k(x;n)
\notag\\[1mm]
&\:\:=
\sum_{k=0}^{n}\frac{(x)_k(1-x)_k(-n)_k}{(2-x)_{k}(1+x)_{k}(2+n)_{k}}
\Big\{H_{k}(x-1)-H_{k}(-x)+H_{k}(1-x)-H_{k}(x)\Big\}.
\end{align*}
Dividing both sides by $1-2x$ and then taking
$(x,n)\to(\frac{1}{2},\infty)$, we have
\begin{align}
\sum_{k=0}^{\infty}\bigg(\frac{1}{4}\bigg)^k\frac{(1)_{k}^2}{(\frac{5}{4})_{k}(\frac{7}{4})_{k}}
\bigg\{\frac{5+6k}{1+2k}\Big[H_{k}^{(2)}+4\Big]-\frac{24+32k}{(1+2k)^3}\bigg\}=24G-24\beta(4).
 \label{eq:wei-bt}
\end{align}
The $(x,n)\to(\frac{1}{2},\infty)$ case of \eqref{eq:wei-at} (cf.
\cite[Example 84]{Chu-b}) reads
\begin{align}
\sum_{k=0}^{\infty}\bigg(\frac{1}{4}\bigg)^k\frac{(1)_{k}^2}{(\frac{5}{4})_{k}(\frac{7}{4})_{k}}\frac{5+6k}{1+2k}=6G.
 \label{eq:wei-ct}
\end{align}
So the combination of \eqref{eq:wei-bt} with \eqref{eq:wei-ct} gives
\eqref{eq:wei-c}.
\end{proof}

Secondly, we plan to prove Theorem \ref{thm-d}.

\begin{proof}[{\bf{Proof of Theorem \ref{thm-d}}}]
Fixing $(a,b,c,d,e)=(1,x,\frac{1}{2},1-x,-n)$ in the hypergeometric
 transformation (cf. \cite[Theorem 14]{Chu-b}):
\begin{align*}
&\sum_{k=0}^{\infty}\frac{(c)_k(e)_k(1+a-b-c)_k(1+a-b-e)_{k}(1+a-c-d)_{k}(1+a-d-e)_{k}}{(1+a-c)_{k}(1+a-e)_{k}}
\notag\\[1mm]
&\quad\times\frac{(1+a-b-d)_{2k}\,\omega_k(a,b,c,d,e)}{(1+a-b)_{2k}(1+a-d)_{2k}(1+2a-b-c-d-e)_{2k}}
\end{align*}
\begin{align*}
&=\sum_{k=0}^{\infty}(a+2k)\frac{(b)_k(c)_k(d)_k(e)_k}{(1+a-b)_{k}(1+a-c)_{k}(1+a-d)_{k}(1+a-e)_{k}},
\end{align*}
where
\begin{align*}
\omega_k(a,b,c,d,e)&=\frac{(1+2a-b-c-d+3k)(a-e+k)}{1+2a-b-c-d-e+2k}
\\[1mm]
&\quad+\frac{(e+k)(1+a-b-c+k)}{(1+a-b+2k)(1+a-d+2k)}
\\[1mm]
&\quad\times\frac{(1+a-c-d+k)(1+a-b-d+2k)(2+2a-b-d-e+3k)}{(1+2a-b-c-d-e+2k)(2+2a-b-c-d-e+2k)},
\end{align*}
 it is ordinary to find that
\begin{align}
&\sum_{k=0}^{n}\frac{(\frac{1}{2}+x)_k(\frac{3}{2}-x)_k(\frac{1}{2})_k(1)_{2k}}{(2-x)_{2k}(1+x)_{2k}(\frac{3}{2})_{k}}
\frac{(1+x+n)_{k}(2-x+n)_{k}(-n)_k}{(2+n)_{k}(\frac{3}{2}+n)_{2k}}E_k(x;n)
\notag\\[1mm]
&\:\:
=\sum_{k=0}^{n}\frac{(x)_k(1-x)_k(-n)_k}{(2-x)_{k}(1+x)_{k}(2+n)_{k}},
\label{eq:wei-dt}
\end{align}
where
\begin{align*}
E_k(x;n) &=\frac{3(1+2k)(1+k+n)}{3+4k+2n}
\notag\\[1mm]
&\quad
+\frac{(1+2x+2k)(3-2x+2k)(1+2k)(3+3k+n)(k-n)}{(1+x+2k)(2-x+2k)(3+4k+2n)(5+4k+2n)}.
\end{align*}
Apply the operator $\mathcal{D}_{x}$ on \eqref{eq:wei-dt} to get
\begin{align*}
&\sum_{k=0}^{n}\frac{(\frac{1}{2}+x)_k(\frac{3}{2}-x)_k(\frac{1}{2})_k(1)_{2k}}{(2-x)_{2k}(1+x)_{2k}(\frac{3}{2})_{k}}
\frac{(1+x+n)_{k}(2-x+n)_{k}(-n)_k}{(2+n)_{k}(\frac{3}{2}+n)_{2k}}E_k(x;n)
\notag\\[1mm]
&\quad\times\Big\{H_{k}(x-\tfrac{1}{2})-H_{k}(\tfrac{1}{2}-x)+H_{2k}(1-x)-H_{2k}(x)+H_{k}(x+n)-H_{k}(1-x+n)\Big\}
\notag\\[1mm]
&\:+\sum_{k=0}^{n}\frac{(\frac{1}{2}+x)_k(\frac{3}{2}-x)_k(\frac{1}{2})_k(1)_{2k}}{(2-x)_{2k}(1+x)_{2k}(\frac{3}{2})_{k}}
\frac{(1+x+n)_{k}(2-x+n)_{k}(-n)_k}{(2+n)_{k}(\frac{3}{2}+n)_{2k}}
\mathcal{D}_{x}E_k(x;n)
\notag\\[1mm]
&\:\:=
\sum_{k=0}^{n}\frac{(x)_k(1-x)_k(-n)_k}{(2-x)_{k}(1+x)_{k}(2+n)_{k}}
\Big\{H_{k}(x-1)-H_{k}(-x)+H_{k}(1-x)-H_{k}(x)\Big\}.
\end{align*}
Dividing both sides by $1-2x$ and then taking
$(x,n)\to(\frac{1}{2},\infty)$, we deduce
\begin{align}
&\sum_{k=0}^{\infty}\bigg(\frac{-1}{4}\bigg)^k\frac{(1)_k^3(\frac{1}{2})_k}{(\frac{5}{4})_k^2(\frac{7}{4})_k^2}
\Big\{(19+56k+40k^2)\big[H_{k}^{(2)}+H_{1+2k}^{(2)}-4H_{3+4k}^{(2)}+4\big]+4\Big\}
\notag\\[1mm]
&\:\: =72G-72\beta(4).
 \label{eq:wei-et}
\end{align}
The $(x,n)\to(\frac{1}{2},\infty)$ case of \eqref{eq:wei-dt} (cf.
\cite[Example 85]{Chu-b}) produces
\begin{align}
\sum_{k=0}^{\infty}\bigg(\frac{-1}{4}\bigg)^k\frac{(1)_k^3(\frac{1}{2})_k}{(\frac{5}{4})_k^2(\frac{7}{4})_k^2}
(19+56k+40k^2)=18G.
 \label{eq:wei-ft}
\end{align}
Hence, the combination of \eqref{eq:wei-et} and \eqref{eq:wei-ft}
leads us to \eqref{eq:wei-d}.
\end{proof}

Thirdly, we start to prove Theorem \ref{thm-e}.

\begin{proof}[{\bf{Proof of Theorem \ref{thm-e}}}]
For achieving the purpose, we need the hypergeometric transformation
 (cf. \cite[Theorem 24]{Chu-b}):
\begin{align}
&\sum_{k=0}^{\infty}\frac{(c)_k(1+a-b-d)_k(1+a-b-e)_k(d)_{2k}(e)_{2k}(1+a-b-c)_{2k}}
{(d+e-a)_{k}(1+a-d)_{k}(1+a-e)_{k}(1+2a-b-c-d-e)_{k}(1+a-c)_{2k}}
\notag\\[1mm]
&\quad\times\frac{(-1)^k}{(1+a-b)_{3k}}\gamma_k(a,b,c,d,e)
\notag\\[1mm]
&=\sum_{k=0}^{\infty}(a+2k)\frac{(b)_k(c)_k(d)_k(e)_k}{(1+a-b)_{k}(1+a-c)_{k}(1+a-d)_{k}(1+a-e)_{k}},
\label{equation-c}
\end{align}
where
\begin{align*}
\gamma_k(a,b,c,d,e)&=\frac{(1+2a-b-c-d+3k)(a-e+k)}{1+2a-b-c-d-e+k}
\\[1mm]
&\quad+\frac{(e+2k)(1+a-b-c+2k)(1+a-b-d+k)}{(1+a-b+3k)(1+2a-b-c-d-e+k)}
\\[1mm]
&\quad+\frac{(c+k)(d+2k)(e+2k)(1+a-b-c+2k)}{(d+e-a+k)(1+a-c+2k)(1+2a-b-c-d-e+k)}
\\[1mm]
&\qquad\times\frac{(1+a-b-d+k)(1+a-b-e+k)}{(1+a-b+3k)(2+a-b+3k)}.
\end{align*}
Setting $(a,b,c,d,e)=(\frac{3}{2},1,-n,x,2-x)$ in
\eqref{equation-c}, we have
\begin{align}
&\sum_{k=0}^{n}\bigg(\frac{-1}{27}\bigg)^k\frac{(-\frac{1}{2}+x)_k(\frac{3}{2}-x)_k(x)_{2k}(2-x)_{2k}}
{(\frac{5}{2}-x)_{k}(\frac{1}{2}+x)_{k}(\frac{1}{2})_{k}^2(\frac{5}{6})_{k}(\frac{7}{6})_{k}}
\frac{(-n)_k(\frac{3}{2}+n)_{2k}}{(1+n)_k(\frac{5}{2}+n)_{2k}}F_k(x;n)
\notag\\[1mm]
&\:\:
=\sum_{k=0}^{n}(3+4k)\frac{(x)_k(2-x)_k(1)_k(-n)_k}{(\frac{5}{2}-x)_{k}(\frac{1}{2}+x)_{k}(\frac{3}{2})_{k}(\frac{5}{2}+n)_{k}},
\label{eq:wei-gt}
\end{align}
where
\begin{align*}
F_k(x;n)
&=\frac{(2x-1+2k)(3-x+3k+n)}{1+k+n}+\frac{(2-x+2k)(3-2x+2k)(3+4k+2n)}{3(1+2k)(1+k+n)}
\\&\quad+\frac{4(2-x+2k)(3-2x+2k)(2x-1+2k)(k-n)(3+4k+2n)}{3(1+2k)^2(5+6k)(1+k+n)(5+4k+2n)}.
\end{align*}
Employ the operator $\mathcal{D}_{x}$ on \eqref{eq:wei-gt} to obtain
\begin{align*}
&\sum_{k=0}^{n}\bigg(\frac{-1}{27}\bigg)^k\frac{(-\frac{1}{2}+x)_k(\frac{3}{2}-x)_k(x)_{2k}(2-x)_{2k}}
{(\frac{5}{2}-x)_{k}(\frac{1}{2}+x)_{k}(\frac{1}{2})_{k}^2(\frac{5}{6})_{k}(\frac{7}{6})_{k}}
\frac{(-n)_k(\frac{3}{2}+n)_{2k}}{(1+n)_k(\frac{5}{2}+n)_{2k}}F_k(x;n)
\notag\\[1mm]
&\quad\times\Big\{H_{2k}(x-1)-H_{2k}(1-x)+H_{k}(x-\tfrac{3}{2})-H_{k}(\tfrac{1}{2}-x)+H_{k}(\tfrac{3}{2}-x)-H_{k}(x-\tfrac{1}{2})\Big\}
\notag\\[1mm]
&\:+\sum_{k=0}^{n}\bigg(\frac{-1}{27}\bigg)^k\frac{(-\frac{1}{2}+x)_k(\frac{3}{2}-x)_k(x)_{2k}(2-x)_{2k}}
{(\frac{5}{2}-x)_{k}(\frac{1}{2}+x)_{k}(\frac{1}{2})_{k}^2(\frac{5}{6})_{k}(\frac{7}{6})_{k}}
\frac{(-n)_k(\frac{3}{2}+n)_{2k}}{(1+n)_k(\frac{5}{2}+n)_{2k}}
\mathcal{D}_{x}F_k(x;n)
\notag\\[1mm]
&\:\:=
\sum_{k=0}^{n}(3+4k)\frac{(x)_k(2-x)_k(1)_k(-n)_k}{(\frac{5}{2}-x)_{k}(\frac{1}{2}+x)_{k}(\frac{3}{2})_{k}(\frac{5}{2}+n)_{k}}
\notag\\[1mm]
&\qquad\times
\Big\{H_{k}(x-1)-H_{k}(1-x)+H_{k}(\tfrac{3}{2}-x)-H_{k}(x-\tfrac{1}{2})\Big\}.
\end{align*}
Dividing both sides by $2-2x$ and then taking $(x,n)\to(1,\infty)$,
it is clear to realize that
\begin{align}
&\sum_{k=0}^{\infty}\bigg(\frac{16}{27}\bigg)^k\frac{(1)_k^2}{(\frac{7}{6})_k(\frac{11}{6})_k}
\bigg\{\frac{21+22k}{1+2k}\Big[H_{1+2k}^{(2)}+4\Big]-\frac{145+174k}{(1+2k)^3}\bigg\}
\notag\\
&\:=
-60\sum_{k=0}^{\infty}\frac{3+4k}{(-1)^k}\frac{(1)_k^3}{(\frac{3}{2})_k^3}
\bigg\{H_{1+2k}^{(2)}-\frac{1}{2}H_{k}^{(2)}-1\bigg\}.
 \label{eq:wei-ht}
\end{align}
The $(x,n)\to(1,\infty)$ case of \eqref{eq:wei-gt} (cf.
\cite[Example 50]{Chu-b}) becomes
\begin{align}
\sum_{k=0}^{\infty}\bigg(\frac{16}{27}\bigg)^k\frac{(1)_k^2}{(\frac{7}{6})_k(\frac{11}{6})_k}
\frac{21+22k}{1+2k}=30G.
 \label{eq:wei-it}
\end{align}
According to \eqref{Guillera-b}, \eqref{Sun-b}, \eqref{eq:wei-ht},
and \eqref{eq:wei-it}, we catch hold of \eqref{eq:wei-e}.
\end{proof}

Finally, we shall prove Theorem \ref{thm-f}.

\begin{proof}[{\bf{Proof of Theorem \ref{thm-f}}}]
Letting $(a,b,c,d,e)=(\frac{3}{2},1,x,2-x,-n)$ in the hypergeometric
transformation (cf. \cite[Theorem 18]{Chu-b}):
\begin{align}
&\sum_{k=0}^{\infty}\frac{(c)_k(d)_k(1+a-b-e)_k(1+a-c-d)_k(e)_{2k}(1+a-b-c)_{2k}(1+a-b-d)_{2k}}
{(1+a-e)_{k}(1+a-c)_{2k}(1+a-d)_{2k}(1+2a-b-c-d-e)_{2k}}
\notag\\[1mm]
&\quad\times\frac{(-1)^k}{(1+a-b)_{3k}}\theta_k(a,b,c,d,e)
\notag\\[1mm]
&=\sum_{k=0}^{\infty}(a+2k)\frac{(b)_k(c)_k(d)_k(e)_k}{(1+a-b)_{k}(1+a-c)_{k}(1+a-d)_{k}(1+a-e)_{k}},
\label{equation-d}
\end{align}
where
\begin{align*}
\theta_k(a,b,c,d,e)&=\frac{(1+2a-b-c-d+4k)(a-e+k)}{1+2a-b-c-d-e+2k}
\\[1mm]
&\quad+\frac{(e+2k)(1+a-b-c+2k)(1+a-b-d+2k)}{(1+a-b+3k)(1+a-d+2k)}
\\[1mm]
&\qquad\times\frac{(1+a-c-d+k)(2+2a-b-d-e+3k)}{(1+2a-b-c-d-e+2k)(2+2a-b-c-d-e+2k)}
\\[1mm]
&\quad+\frac{(c+k)(e+2k)(1+a-b-c+2k)(1+a-b-e+k)}{(1+a-c+2k)(1+a-d+2k)(1+a-b+3k)(2+a-b+3k)}
\\[1mm]
&\qquad\times\frac{(1+a-c-d+k)(1+a-b-d+2k)(2+a-b-d+2k)}{(1+2a-b-c-d-e+2k)(2+2a-b-c-d-e+2k)},
\end{align*}
it is obvious to provide that
\begin{align}
&\sum_{k=0}^{n}\bigg(\frac{-1}{27}\bigg)^k\frac{(x)_{k}(2-x)_{k}(-\frac{1}{2}+x)_{2k}(\frac{3}{2}-x)_{2k}}
{(\frac{5}{2}-x)_{2k}(\frac{1}{2}+x)_{2k}(\frac{5}{6})_{k}(\frac{7}{6})_{k}}
\frac{(\frac{3}{2}+n)_{k}(-n)_{2k}}{(\frac{5}{2}+n)_{k}(1+n)_{2k}}G_k(x;n)
\notag\\[1mm]
&\:\:
=\sum_{k=0}^{n}(3+4k)\frac{(x)_k(2-x)_k(1)_k(-n)_k}{(\frac{5}{2}-x)_{k}(\frac{1}{2}+x)_{k}(\frac{3}{2})_{k}(\frac{5}{2}+n)_{k}},
\label{eq:wei-jt}
\end{align}
where
\begin{align*}
G_k(x;n)
&=\frac{(1+4k)(3+2k+2n)}{1+2k+n}
\\&\quad+
\frac{(2x-1+4k)(3-2x+4k)(2k-n)(2+x+3k+n)}{3(1+2x+4k)(1+2k+n)(2+2k+n)}
\\&\quad+\frac{(x+k)(2x-1+4k)(3-2x+4k)(2k-n)(3+2k+2n)}{3(5-2x+4k)(5+6k)(1+2k+n)(2+2k+n)}.
\end{align*}
Apply the operator $\mathcal{D}_{x}$ on \eqref{eq:wei-jt} to get
\begin{align*}
&\sum_{k=0}^{n}\bigg(\frac{-1}{27}\bigg)^k\frac{(x)_{k}(2-x)_{k}(-\frac{1}{2}+x)_{2k}(\frac{3}{2}-x)_{2k}}
{(\frac{5}{2}-x)_{2k}(\frac{1}{2}+x)_{2k}(\frac{5}{6})_{k}(\frac{7}{6})_{k}}
\frac{(\frac{3}{2}+n)_{k}(-n)_{2k}}{(\frac{5}{2}+n)_{k}(1+n)_{2k}}G_k(x;n)
\notag\\[1mm]
&\quad\times\Big\{H_{k}(x-1)-H_{k}(1-x)+H_{2k}(x-\tfrac{3}{2})-H_{2k}(\tfrac{1}{2}-x)+H_{2k}(\tfrac{3}{2}-x)-H_{2k}(x-\tfrac{1}{2})\Big\}
\notag\\[1mm]
&\:+\sum_{k=0}^{n}\bigg(\frac{-1}{27}\bigg)^k\frac{(x)_{k}(2-x)_{k}(-\frac{1}{2}+x)_{2k}(\frac{3}{2}-x)_{2k}}
{(\frac{5}{2}-x)_{2k}(\frac{1}{2}+x)_{2k}(\frac{5}{6})_{k}(\frac{7}{6})_{k}}
\frac{(\frac{3}{2}+n)_{k}(-n)_{2k}}{(\frac{5}{2}+n)_{k}(1+n)_{2k}}\mathcal{D}_{x}G_k(x;n)
\notag\\[1mm]
&\:\:=
\sum_{k=0}^{n}(3+4k)\frac{(x)_k(2-x)_k(1)_k(-n)_k}{(\frac{5}{2}-x)_{k}(\frac{1}{2}+x)_{k}(\frac{3}{2})_{k}(\frac{5}{2}+n)_{k}}
\notag\\[1mm]
&\qquad\times
\Big\{H_{k}(x-1)-H_{k}(1-x)+H_{k}(\tfrac{3}{2}-x)-H_{k}(x-\tfrac{1}{2})\Big\}.
\end{align*}
Dividing both sides by $2-2x$ and then taking $(x,n)\to(1,\infty)$,
we discover
\begin{align}
&\sum_{k=0}^{\infty}\bigg(\frac{-1}{27}\bigg)^k\frac{(1)_k^2}{(\frac{7}{6})_k(\frac{11}{6})_k}
\bigg\{\frac{83+192k+112k^2}{64(1+4k)(3+4k)}\Big[H_{k}^{(2)}+4\Big]
\notag\\
&\qquad\:\:
-\frac{(5+6k)(10+39k+48k^2+16k^3)}{(1+4k)^3(3+4k)^3}\bigg\}
\notag\\
&\:=
-\frac{15}{16}\sum_{k=0}^{\infty}\frac{3+4k}{(-1)^k}\frac{(1)_k^3}{(\frac{3}{2})_k^3}
\bigg\{H_{1+2k}^{(2)}-\frac{1}{2}H_{k}^{(2)}-1\bigg\}.
 \label{eq:wei-lt}
\end{align}
The $(x,n)\to(1,\infty)$ case of \eqref{eq:wei-jt} (cf.
\cite[Example 29]{Chu-b}) engenders
\begin{align}
\sum_{k=0}^{\infty}\bigg(\frac{-1}{27}\bigg)^k\frac{(1)_k^2}{(\frac{7}{6})_k(\frac{11}{6})_k}
\frac{83+192k+112k^2}{(1+4k)(3+4k)}=30G.
 \label{eq:wei-mt}
\end{align}
In terms of  \eqref{Guillera-b}, \eqref{Sun-b}, \eqref{eq:wei-lt},
and \eqref{eq:wei-mt}, we are led to \eqref{eq:wei-f}.
\end{proof}


\end{document}